\documentclass{article}
\usepackage[leqno]{amsmath}
\usepackage{amsmath,amsfonts,amsthm,amssymb}
\usepackage{epsfig}
\usepackage{graphicx}
\usepackage{subfig}

\usepackage{cite}
\usepackage{afterpage,hyperref}
\hypersetup{
  colorlinks,
  citecolor=blue!50!black,
  linkcolor=red,
  urlcolor=blue!50!black}
\usepackage{tikz}
\usetikzlibrary{backgrounds,calc,positioning}
\usepackage{epsfig}
\usepackage{color}
\usepackage{amsmath}
\usepackage{amssymb}
\newtheorem{lemma}{Lemma}

\newtheorem{theorem}{Theorem}
\newtheorem{corollary}{Corollary}

\theoremstyle{definition}

\newtheorem{conjecture}{Conjecture}

\usetikzlibrary{decorations.markings} 

\tikzstyle{graphnode}=[draw,shape=circle,draw=black,minimum size=0.5pt,inner sep=1.5pt]

\title{$(2,3)$-bipartite graphs are strongly 6-edge-choosable}
\author{Petru Valicov~\thanks{Aix Marseille Univ, Universit\'e de Toulon, CNRS, LIS, Marseille, France.\newline
\indent Email: petru.valicov@lis-lab.fr}}

\begin{document}
\maketitle

\begin{abstract}
Kang and Park recently showed that every cubic (loopless) multigraph is incidence 6-choosable [On incidence choosability of cubic graphs. \emph{arXiv}, April 2018]. Equivalently, every bipartite graph obtained by subdividing once every edge of a cubic graph, is strongly 6-edge-choosable. The aim of this note is to give a shorter proof of their result by looking at the strong edge-coloring formulation of the problem.
\end{abstract}

In this note all graphs are finite without loops. For a graph $G$, we denote by $\Delta(G)$ (or simply $\Delta$ when no ambiguity is possible) the \textit{maximum degree} of $G$. A \textit{strong edge-coloring} of $G$ is a proper edge-coloring of $G$ such that every two edges joined by another edge are colored differently. The least number of colors in a strong edge-coloring of $G$ is referred to as the \textit{strong chromatic index}, denoted $\chi'_s(G)$.

\begin{conjecture}[Brualdi and Quinn Massey~\cite{BQM93}]
\label{conj:brualdi_bipartite}
For every bipartite graph $G$ with bipartition $A$ and $B$, we have $\chi'_s(G)\leq \Delta(A)\Delta(B)$.
\end{conjecture}

Define a \emph{$(d_A, d_B)$-bipartite graph} to be a bipartite graph with parts $A$ and $B$ such that $\Delta(A)\leq d_A$ and $\Delta(B)\leq d_B$. Conjecture~\ref{conj:brualdi_bipartite} was shown to be true for $(2,\Delta)$-bipartite graphs~\cite{N08} and for $(3,\Delta)$-bipartite graphs~\cite{BLV16,HYZ17}.

As for the classical proper coloring, one can define the list-version of the problem.
For an edge $uv$, let $L(uv)$ denote the list of available colors. A list coloring is a choice function that maps every edge $uv$ to a color in the list $L(uv)$. Denote by $ch'_s(G)$ the strong-edge choosability of $G$ i.e. the smallest number $k$ such that $G$ admits a list strong edge-coloring with $|L(uv)| \leq k, \forall uv\in E(G)$.

Note that no counterexample for list-coloring version of Conjecture~\ref{conj:brualdi_bipartite} is known. The aim of this note is to show the following result (already proved by Kang and Park in~\cite{KP18} using different techniques):

\begin{theorem}
\label{thm:main}
For every $(2, 3)$-bipartite graph $G$, we have $ch'_s(G)\leq 6$.
\end{theorem}

Observe that the bound is tight, the complete bipartite graph $K_{2,3}$ being an example. 

Let us now recall the definition of \emph{incidence coloring} which is a special case of strong edge-coloring. Let $G$ be an undirected graph and let $v$ be a vertex of $G$ incident to edge $e$. We call $(v,e)$ an incidence. Two distinct incidences $(v,e)$ and $(w,f)$ are adjacent if one of the following holds:
\begin{itemize}
\item $v = w$ or
\item $e = f$ or
\item $vw \in \{e, f\}$
\end{itemize}
An \emph{incidence coloring} of $G$ assigns a color to each incidence of $G$ such that adjacent incidences get distinct colors. The smallest number of colors required for such a coloring is called the \emph{incidence chromatic number} of $G$, and is denoted by $\chi_i(G)$. The notion of incidence coloring was introduced in 1993 by Brualdi and Quinn Massey~\cite{BQM93} and is a particular case of strong edge-coloring: given a graph $G$ with maximum degree $\Delta$, an incidence coloring of $G$ corresponds to a strong edge-coloring of the $(2,\Delta)$-bipartite graph obtained from $G$ by subdividing each edge of $G$ exactly once. Therefore, Theorem~\ref{thm:main} implies the following:

\begin{corollary}
\label{cor:main}
Let $G$ be a multigraph with $\Delta\leq 3$. Then $G$ admits an incidence coloring with at most 6 colors.
\end{corollary}

\begin{figure}[!ht]
\centering
\subfloat[Initial configuration]
{
\label{subfig:config1}
\scalebox{0.85}{
\begin{tikzpicture}[join=bevel,inner sep=0.4mm]
  \foreach \i / \l in {0/u,2/v,4/w,6/x,8/y}{
  	\node[circle,draw] (\l) at (\i,0) {$\l$};
  }
  \foreach \i in {1,3,5,7}{
  	\node at (\i,-0.3) {$\geq 5$};
  }
  \node[circle,draw] (z) at (2,2) {$z$};
  \node at (1.6,1) {$\geq 3$};
  \node[circle,draw] (t) at (6,2) {$t$};
  \node at (5.6,1) {$\geq 3$};
  \draw[-] (u) -- (v) -- (w) -- (x) -- (y);
  \draw[-] (v) -- (z);
  \draw[-] (x) -- (t);
\end{tikzpicture}
}
}
\subfloat[The configuration after coloring $uv,vz,xt,xy$]
{
\label{subfig:config1Precol}
\scalebox{0.85}{
\begin{tikzpicture}[join=bevel,inner sep=0.4mm]
  \foreach \i / \l in {0/u,2/v,4/w,6/x,8/y}{
  	\node[circle,draw] (\l) at (\i,0) {$\l$};
  }
  \node at (3,-0.3) {$\geq 3$};
  \node at (5,-0.3) {$\geq 2$};
  \node[circle,draw] (z) at (2,2) {$z$};
  \node[circle,draw] (t) at (6,2) {$t$};
  \draw[-] (u) -- (v) -- (w) -- (x) -- (y);
  \draw[-] (v) -- (z);
  \draw[-] (x) -- (t);
  
\end{tikzpicture}
}
}
 \caption{The graph from Lemma~\ref{lem:config1}}
\label{fig:config1}
\end{figure}
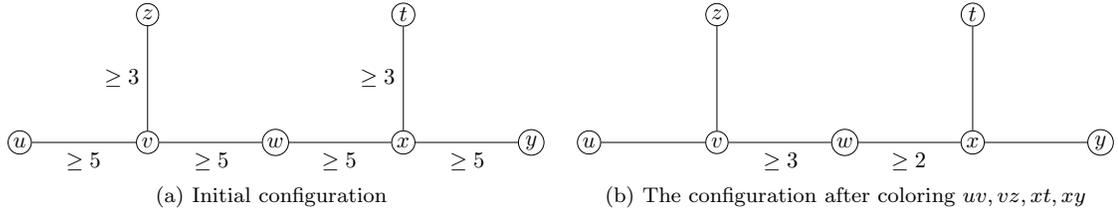

To show Theorem~\ref{thm:main}, we need two additional results:

\begin{lemma}
\label{lem:config1}
Let $G$ be a $(2,3)$-bipartite graph with 6 vertices formed from a path on 5 vertices $uvwxy$ and with two additional vertices $z,t$ adjacent to $v$ and $x$, respectively (see Figure~\ref{subfig:config1}). Let $|L(uv)|\geq 5,|L(vw)|\geq 5, |L(wx)|\geq 5,|L(xy)|\geq 5,|L(vz)|\geq 3, |L(xt)|\geq 3$. Then it is always possible to assign colors to edges $uv,vz,xy,xt$ from their respective lists, such that the obtained precoloring is a valid strong edge-coloring and moreover we have $|L(vw)|\geq 3, |L(wx)|\geq 2$ (see Figure~\ref{subfig:config1Precol}).
\end{lemma}

\begin{proof}
We distinguish two cases:
\begin{itemize}
\item Suppose $L(vz)\cap L(xt)\neq \emptyset$. Assign then $vz$ and $xt$ a color $\alpha\in L(vz)\cap L(xt)$. Now the lists of the remaining uncolored edges have the following sizes: $|L(uv)|\geq 4,|L(vw)|\geq 4, |L(wx)|\geq 4,|L(xy)|\geq 4$. By looking at $L(uv)\cap L(xy)$ (whether the intersection is empty or not), one can easily assign colors to $uv,xy$ such that at the end we have $|L(vw)|\geq 3, |L(wx)|\geq 2$.

\item Suppose $L(vz)\cap L(xt) = \emptyset$ and thus $|L(vz)\cup L(xt)|\geq 6$. We denote $\{e_1,e_2\}=\{vz,xt\}$. By pigeonhole principle, assign a color to one of the edges of $\{e_1,e_2\}$ (say $e_1$) such that after doing so, we have $|L(vw)|\geq  5$ and $|L(wx)|\geq 4$. Now if $L(uv)\cap L(xy)\neq \emptyset$, then by assigning the same color to both $uv$ and $xy$ and coloring $e_2$, we are done. Thus $L(uv)\cap L(xy) = \emptyset$ and therefore, since $e_1$ was already colored, we have $|L(uv)\cup L(xy)|\geq 9$.
Let $\{e_3,e_4\}=\{uv,xy\}$. By pigeonhole principle one can assign a color $\alpha\in L(uv)\cup L(xy)$ to one of the edges $\{e_3,e_4\}$ (say $e_3$) such that $\alpha\notin L(vw)$. By doing so, we have $|L(vw)|\geq 5$ and $|L(wx)|\geq 3$. Since $|L(e_4)|\geq 4$, we assign a color $\beta\in L(e_4)$ such that we still have $|L(wx)|\geq 3$. Finally, by coloring $e_2$ with one of the remaining colors in $L(e_2)$ we get $|L(vw)|\geq 3, |L(wx)|\geq 2$.
\end{itemize}

\end{proof}

\begin{figure}[!ht]
\centering
\subfloat[Initial configuration]
{
\label{subfig:config2}
\begin{tikzpicture}[join=bevel,state/.style={circle,draw, minimum size=0.5cm},inner sep=0mm]
  \foreach \i / \l in {0/1,2/2,4/3,6/4,8/5}{
  	\node[state] (\l) at (\i,0) {$v_\l$};
  }
  \node at (1,-0.3) {$\geq 3$};
  \node at (3,-0.3) {$\geq 4$};
  \node at (5,-0.3) {$\geq 4$};
  \node at (7,-0.3) {$\geq 3$};
  
  \node[state] (z) at (2,2) {$v'_2$};
  \node at (1.6,1) {$\geq 2$};
  \node[state] (t) at (6,2) {$v'_4$};
  \node at (5.6,1) {$\geq 2$};
  
  \draw[-] (1) -- (2) -- (3) -- (4) -- (5);
  \draw[-] (2) -- (z);
  \draw[-] (4) -- (t);
\end{tikzpicture}
}
\vspace*{1cm}
\subfloat[General configuration]
{
\label{subfig:config2General}

\scalebox{0.85}{
\begin{tikzpicture}[join=bevel,state/.style={circle,draw, minimum size=0.7cm},inner sep=0mm]
  \foreach \i / \l in {0/1,2/2,4/3,6/4,8/5}{
  	\node[state] (\l) at (\i,0) {$v_\l$};
  }
  \node at (1,-0.3) {$\geq 3$};
  \node at (3,-0.3) {$\geq 4$};
  \node at (5,-0.3) {$\geq 5$};
  \node at (7,-0.3) {$\geq 5$};
  
  \node at (10,0) {$\ldots$};
  \foreach \i / \l / \n in {12/3/6,14/2/7,16/1/8}{
  	\node[state] (\n) at (\i,0) {\footnotesize $v_{n-\l}$};
  }
  \node[state] (n) at (18,0) {\footnotesize $v_n$};

  \node at (17,-0.3) {$\geq 3$};
  \node at (15,-0.3) {$\geq 4$};
  \node at (13,-0.3) {$\geq 5$};
  \node at (11,-0.3) {$\geq 5$};  
  \node at (9,-0.3) {$\geq 5$};
  
  \node[state] (v'2) at (2,2) {\footnotesize $v'_2$};
  \node at (1.6,1) {$\geq 2$};
  \node[state] (v'4) at (6,2) {\footnotesize $v'_4$};
  \node at (5.6,1) {$\geq 3$};
  
  \node[state] (v'n-3) at (12,2) {\footnotesize $v'_{n-3}$};
  \node at (15.6,1) {$\geq 2$};
  \node[state] (v'n-1) at (16,2) {\footnotesize $v'_{n-1}$};
  \node at (11.6,1) {$\geq 3$};

  \draw[-] (1) -- (2) -- (3) -- (4) -- (5);
  \draw[-] (5) -- (9.5,0) (10.5,0) -- (6);
  \draw[-] (6) -- (7) -- (8) -- (n);
  \draw[-] (2) -- (z) (4) -- (t);
  \draw[-] (6) -- (v'n-3) (8) -- (v'n-1);
\end{tikzpicture}
}
}
 \caption{The graph from Lemma~\ref{lem:config2}}
\label{fig:config2}
\end{figure}
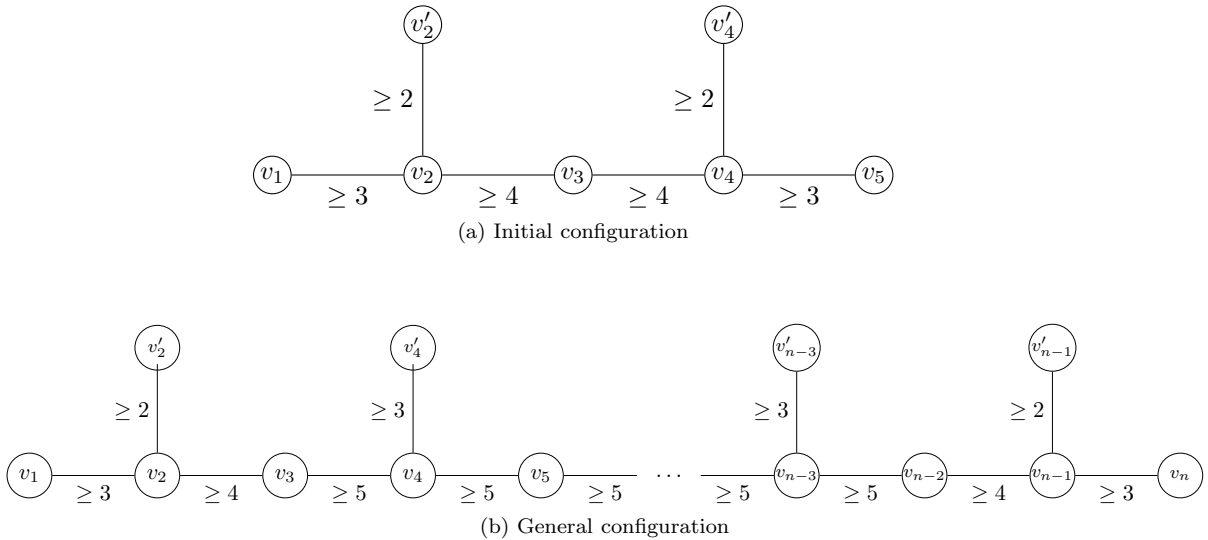

\begin{lemma}
\label{lem:config2}
Let $G$ be a $(2,3)$-bipartite graph formed from a path $v_1v_2\ldots v_n$, with $n\geq 5$ being odd, and for every even $i \in \{2,4,\ldots,n-1\}$, let $v'_i$ be the third vertex adjacent to $v_i$ (see Figure~\ref{fig:config2}). Let $|L(v_1v_2)|\geq 3$, $|L(v_2v'_2)|\geq 2$, $|L(v_2v_3)|\geq 4$, $|L(v_{n-2}v_{n-1})|\geq 4$, $|L(v_{n-1}v'_{n-1})|\geq 2$,  $L(v_{n-1}v_n)|\geq 3$. For $n \geq 7$, let $|L(v_iv_{i+1})|\geq 5$, $\forall i \in \{3,\ldots,n-3\}$ and for every even $j \in \{4,6,\ldots,n-3\}$, let $|L(v_jv'_j)|\geq 3$. Then it is always possible to strong list edge-color $G$.
\end{lemma}

\begin{proof}
We induct on $n$. For $n=5$, the configuration is depicted in Figure~\ref{subfig:config2}. We distinguish the following cases, where in each case we suppose that none of the
previous ones apply:

\begin{itemize}
\item Suppose there exists a color $\alpha\in L(v_2v'_2)\cap L(v_4v'_4)$. Color edges $v_2v'_2, v_4v'_4$ with $\alpha$. Now we have, $|L(v_1v_2)|\geq 2,|L(v_4v_5)|\geq 2,|L(v_2v_3)|\geq 3,|L(v_3v_4)|\geq 3$. By considering $L(v_1v_2)\cap L(v_4v_5)$, one can easily see that the coloring can be extended. Thus, we have that $L(v_2v'_2)\cap L(v_4v'_4)=\emptyset$.

\item Suppose $L(v_1v_2)\cap L(v_4v'_4)\neq \emptyset$. Assign $v_1v_2$ and $v_4v'_4$ the same color $\alpha\in L(v_1v_2)\cap L(v_4v'_4)$. Now we have the following sizes of lists, $|L(v_2v_3)|\geq 3$, $|L(v_3v_4)|\geq 3$, $|L(v_4v_5)|\geq 2$, $|L(v_2v'_2)|\geq 2$. As in the previous case, by considering $L(v_2v'_2)\cap L(v_4v_5)$, one can easily see that the coloring can be extended to $G$. Therefore, $L(v_1v_2)\cap L(v_4v'_4) = \emptyset$ and symmetrically $L(v_2v'_2)\cap L(v_4v_5) = \emptyset$.

\item Suppose there exists a color $\alpha\in L(v_1v_2)\cap L(v_4v_5)$. Assign $\alpha$ to edges $v_1v_2,v_4v_5$. Now, by the previous cases, we have $|L(v_2v'_2)|\geq 2$, $|L(v_4v'_4)|\geq 2$, $|L(v_2v_3)|\geq 3$, $|L(v_3v_4)|\geq 3$. Again, by considering $L(v_2v'_2)\cap L(v_4v'_4)$, one can easily see that the coloring can be extended to $G$.
\end{itemize}

We have $L(v_1v_2)\cap L(v_4v_5) = \emptyset$ and $L(v_1v_2)\cap L(v_4v'_4) = \emptyset$. Hence $|L(v_1v_2)\cup L(v_4v_5)| \geq 6$ and $|L(v_1v_2)\cap L(v_4v'_4)| \geq 5$. Moreover, by symmetry $|L(v_2v'_2)\cup L(v_4v_5)| \geq 5$. Therefore, for every subset $S$ of $\{v_1v_2,v_2v_3,v_3v_4,v_4v_5,v_2v'_2,v_4v'_4\}$ we have $|S|\leq |\bigcup_{e \in S}L(e)|$, and by Hall's Theorem one can list-color all the edges of $G$.

For the inductive step, we have the general configuration depicted in Figure~\ref{subfig:config2General}. Observe that $|L(v_2v_3)|\geq 4$ and $|L(v_4v'_4)|\geq 3$. By pigeonhole principle we assign to $v_2v_3$ a color $\alpha\in L(v_2v_3)$ such that after doing it, $|L(v_4v'_4)|\geq 3$. Then we color greedily $v_2v'_2,v_1v_2$, in this order. We obtain the following modifications on the sizes of the lists: $|L(v_3v_4)|\geq 2$, $|L(v_4v'_4)|\geq 3$, $|L(v_4v_5)|\geq 4$. The lists of other edges are not modified. Therefore, the remaining uncolored edges induce a graph which is smaller and satisfy the induction hypothesis and thus we are done.
\end{proof}

\section*{Proof of Theorem~\ref{thm:main}}

We prove the statement by induction on the order of the graph. Let $G$ be a counterexample to the theorem of the smallest order. By minimality of $G$ it is easy to see $\Delta(A)=2$ and $\Delta(B)=3$ ($A,B$ being the bipartition of $G$).

First we show that $G$ has no 4-cycle. Indeed, if $G$ had a 4-cycle $uvwx$, then since $G$ is a $(2,3)$-bipartite graph, $d(u)\leq 2$ and $d(w)\leq 2$. Let $v'$ (resp. $x'$) be the third neighbor of $v$ (resp. $x$). Consider $G'=G-\{u,v,w,x\}$. By minimality of $G$, $G'$ is 6-choosable. Therefore, for any precoloring of $G'$ in $G$, we have the following lists of available colors for the remaining uncolored edges of $G$ : $|L(vv')|\geq 3, |L(xx')|\geq 3, |L(uv)|\geq 5, |L(vw)|\geq 5, |L(wx)|\geq 5, |L(xv)|\geq 5$. If $|L(vv')\cup L(xx')|\geq 6$, then for every subset $S$ of $\{vv',xx',uv,vw,wx,xv\}$ we have $|S|\leq |\bigcup_{e \in S}L(e)|$, and thus by Hall's theorem, there exists a choice for each remaining uncolored edge of $G$. Hence $L(vv')=L(xx')$. Then assigning $vv',xx'$ the same color from $L(vv')$, one can greedily color the edges of the 4-cycle.

Now, we show that $G$ has no 6-cycle. Let $uvwxyz$ be a 6-cycle in $G$ and $v,x,z$ be the vertices of degree at most 2 and $u',w',y'$ be the neighbors of $u,w,y$ respectively (outside the 6-cycle). Consider the graph $G'=G-\{u,v,w,x,y,z\}$, which by minimality of $G$ is 6-choosable. For any list-coloring of $G'$ in $G$, we have the following lists of available colors for the remaining uncolored edges of $G$: $|L(uu')|\geq 3, |L(ww')|\geq 3, |L(yy')|\geq 3, |L(uv)|\geq 5, |L(vw)|\geq 5, |L(wx)|\geq 5, |L(xy)|\geq 5, |L(yz)|\geq 5, |L(zu)|\geq 5$.
\begin{itemize}
\item Assume $L(uv)\cap L(xy)\neq \emptyset$ and let $\alpha \in L(uv)\cap L(xy)$. Assign to $uv$ and $xy$ color $\alpha$. If there exists a color $\beta\in L(uu')\cap L(yy')$, then assign $\beta$ to $uu',yy'$ and color greedily $ww',vw,wx,yz,zu$ in this order. Hence, $L(uu')\cap L(yy')=\emptyset$ and with the similar reasoning one can easily see that also $L(ww')\cap L(yy')=\emptyset$, $L(ww')\cap L(uu')=\emptyset$. Observe that if there exists $\beta\in L(uu')\cap L(wx)$, then by assigning $\beta$ to $uu'$ and $wx$ and coloring $yy',yz,zu,vw,ww'$ in this order, we are done. Therefore $L(uu')\cap L(wx)=\emptyset$ and similarly one can conclude that $L(yy')\cap L(vwx)=\emptyset$, $L(ww')\cap L(zu)=\emptyset$, $L(ww')\cap L(yz)=\emptyset$. Now, we have that for every subset $S$ of $\{uu',ww',yy',uv,vw,wx,xy,yz,zu\}$ we have $|S|\leq |\bigcup_{e \in S}L(e)|$, and thus by Hall's theorem, there exists a unique choice for each remaining uncolored edge of $G$.
\item We have $L(uv)\cap L(xy) = \emptyset$ and thus $|L(uv)\cup L(xy)|\geq 10$. Symmetrically, $|L(vw)\cup L(yz)|\geq 10$ and  $|L(wx)\cup L(zu)|\geq 10$. Now, if $L(uu')=L(ww')=L(yy')$ then by assigning to $uu',ww',yy'$ the same color, one of the edges of the 6-cycle, say $uv$, will satisfy $|L(uv)|\geq 5$. Color greedily the edges of the 6-cycle in this order : $vw,wx,xy,yz,zu,uv$.
\item Now, observe that for every subset $S$ of $\{uu',ww',yy',uv,vw,wx,xy,yz,zu\}$ we have $|S|\leq |\bigcup_{e \in S}L(e)|$ and thus by Hall's theorem there exists a choice of colors for each of the 9 uncolored edges of $G$.
\end{itemize}

We know now that $G$ has girth (i.e. the length of the shortest cycle) at least 8. Take any cycle in $v_1v_2,\ldots,v_n$ with $n$ even and $n\geq 8$. Assume without loss of generality that $v_1,v_3,v_5,\ldots,v_{n-1}$ have degree 2. Let $G'=G-\{v_1,v_2,\ldots,v_n\}$ and consider any precoloring of $G'$ with respect to $G$. Observe that for every edge $e\in \{v_1,v_2,\ldots,v_n\}$, $|L(e)|\geq 5$ and for every even $i\in\{2,4,\ldots,n\}$,  $|L(v_iv'_i)|\geq 3$.
We provide a list-coloring algorithm in three steps:
\begin{enumerate}
\item By Lemma~\ref{lem:config1}, the edges of the subgraph induced by the vertices $v_1,v_2,v_3,v_4,v_5,v'_2,v'_4$ can be precolored such that edges $v_2v_3$ and $v_3v_4$ are the only edges remaining to color and $|L(v_2v_3)|\geq 3, |L(v_3v_4)|\geq 2$. We do this precoloring.
\item After Step 1, the subgraph induced by the vertices $v_1,v_5,v_{n},v_{n-1},\ldots,v_6,v'_6,v'_8,\ldots,v'_n$ satisfies the conditions of Lemma~\ref{lem:config2} and thus can be entirely list-colored.
\item After Step 2, it remains to color edges $v_2v_3$ and $v_3v_4$. Observe that $|L(v_2v_3)|\geq 2$ and $L(v_3v_4)|\geq 1$. Hence, we color $v_3v_4$, $v_2v_3$ (in this order) and we are done.
\end{enumerate}

%\begin{remark}
%The base cases of the induction in Lemmas~\ref{lem:config1},~\ref{lem:config2} can be easily proved using the Combinatorial Nullstellensatz~\cite{A99}. This can significantly shorten the proofs.
%\end{remark}

\end{document}